\documentclass[11pt,oneside,fleqn]{article}

\usepackage[ansinew]{inputenc}
\usepackage[mathscr]{eucal}
\usepackage{amsmath,amssymb,amsthm}
\usepackage{graphicx}
\usepackage{cite}
\usepackage{hyperref}

\allowdisplaybreaks

\hypersetup{colorlinks, linkcolor=blue, citecolor=blue, urlcolor=blue}
\makeatletter
\long\def\@makecaption#1#2{%
  \vskip\abovecaptionskip\footnotesize
  \sbox\@tempboxa{#1. #2}%
  \ifdim \wd\@tempboxa >\hsize
    #1. #2\par
  \else
    \global \@minipagefalse
    \hb@xt@\hsize{\hfil\box\@tempboxa\hfil}%
  \fi
  \vskip\belowcaptionskip}
\makeatother

\newcommand{\ddd}{\mathrm{d}}
\newcommand{\p}{\partial}
\newcommand{\const}{\mathop{\rm const}\nolimits}

\newcommand{\todo}[1][\null]{\ensuremath{\clubsuit}}

\newcommand{\noprint}[1]{}

{\theoremstyle{definition}

\newtheorem*{remark*}{Remark}
}

\newcommand{\checked}[1][\null]{\ensuremath{\boldsymbol{\surd}}}

\newcommand{\ve}{\varepsilon}

\begin{document}

\par\noindent {\LARGE\bf
Convecting reference frames and \\invariant numerical models
\par}

{\vspace{4mm}\par\noindent {Alexander Bihlo$^\dag\hspace{0.1mm}^\ddag$ and Jean--Christophe Nave$^\ddag$
} \par\vspace{2mm}\par}

{\vspace{2mm}\par\noindent {\it
$^{\dag}$~Centre de recherches math\'{e}matiques, Universit\'{e} de Montr\'{e}al, C.P.\ 6128, succ.\ Centre-ville,\\
$\phantom{^\dag}$~Montr\'{e}al (QC) H3C 3J7, Canada
}}

{\vspace{2mm}\par\noindent {\it
$^{\ddag}$~Department of Mathematics and Statistics, McGill University, 805 Sherbrooke W.,\\
$\phantom{^\dag}$~Montr\'{e}al (QC) H3A 2K6, Canada
}}

{\vspace{2mm}\par\noindent {\it
$\phantom{^\dag}$~\textup{E-mail}:
bihlo@crm.umontreal.ca, jcnave@math.mcgill.ca
}\par}

\vspace{4mm}\par\noindent\hspace*{8mm}\parbox{140mm}{\small
 In the recent paper by Bernardini et al.\ [{\it J.\ Comput.\ Phys.} \textbf{232} (2013), 1--6] the discrepancy in the performance of finite difference and spectral models for simulations of flows with a preferential direction of propagation was studied. In a simplified investigation carried out using the viscous Burgers equation the authors attributed the poorer numerical results of finite difference models to a violation of Galilean invariance in the discretization and propose to carry out the computations in a reference frame moving with the bulk velocity of the flow. Here we further discuss this problem and relate it to known results on invariant discretization schemes. Non-invariant and invariant finite difference discretizations of Burgers equation are proposed and compared with the discretization using the remedy proposed by Bernardini et al..
}\par\vspace{2mm}

\section{Introduction}

In the recent paper~\cite{bern13Ay} a possible remedy was discussed to improve the poor numerical behavior of finite difference simulations of turbulent flows with a preferential propagation direction. It was shown that the violation of Galilean invariance of the finite difference scheme is the most likely explanation why it is necessary to use a significantly larger number of grid points in finite difference calculations than in spectral methods to achieve comparably accurate numerical results. The recommendation given in~\cite{bern13Ay} is to carry out the finite difference computations in a reference frame that moves with the constant stream-wise bulk velocity in the flow direction. It was then shown for the example of Burgers equation that the finite difference model may yield similar numerical results as spectral discretizations with approximately the same number of grid points.

In the present paper we further discuss this problem and the remedy proposed in~\cite{bern13Ay}. In fact, the problem found and analyzed in~\cite{bern13Ay} has been investigated quite intensively in the field of group analysis of differential and difference equations, see e.g.~\cite{bihl12Cy,bihl12By,chha11Ay,doro11Ay,kim08Ay,levi06Ay,rebe11Ay} and references therein for some of the most recent results. In particular, it was established by Dorodnitsyn and collaborators~\cite{budd01Ay,doro11Ay,doro03Ay} that it is not possible to maintain the Galilean invariance of partial differential equations in a finite difference model when the mesh does not move in the course of the numerical integration. This result qualitatively explains why the method proposed in~\cite{bern13Ay} may work from the geometrical point of view.

The violation of Galilean invariance of stationary discretizations can be readily shown by applying a Galilean boost, which in the one-dimensional case is
\begin{equation}\label{eq:GalileanTransformation}
 (\tilde t,\tilde x,\tilde u)=(t,x+\ve t,u+\ve),
\end{equation}
where $\ve\in\mathbb{R}$, to the defining equation of the grid,
$
 x^{n+1}_i-x^n_i=0.
$
Here and in the following, an upper index indicates the time level and a lower index the spatial grid point. The action of the Galilean transformation~\eqref{eq:GalileanTransformation} on this grid equation yields
$
 \tilde x^{n+1}_i-\tilde x^n_i=x^{n+1}_i-x^n_i+\ve (t^{n+1}-t^n),
$
which clearly fails to be invariant for $\ve\ne0$. Here we assumed that all the grid points are defined on the same time layer, i.e.\ $t^n_{i+1}=t^n_i=t^n$. It can be checked that this assumption does not violate the invariance of most of the equations of hydrodynamics, see also~\cite{doro11Ay} for more details.

\begin{figure}[!ht]
 \centering
 \includegraphics[scale=0.40]{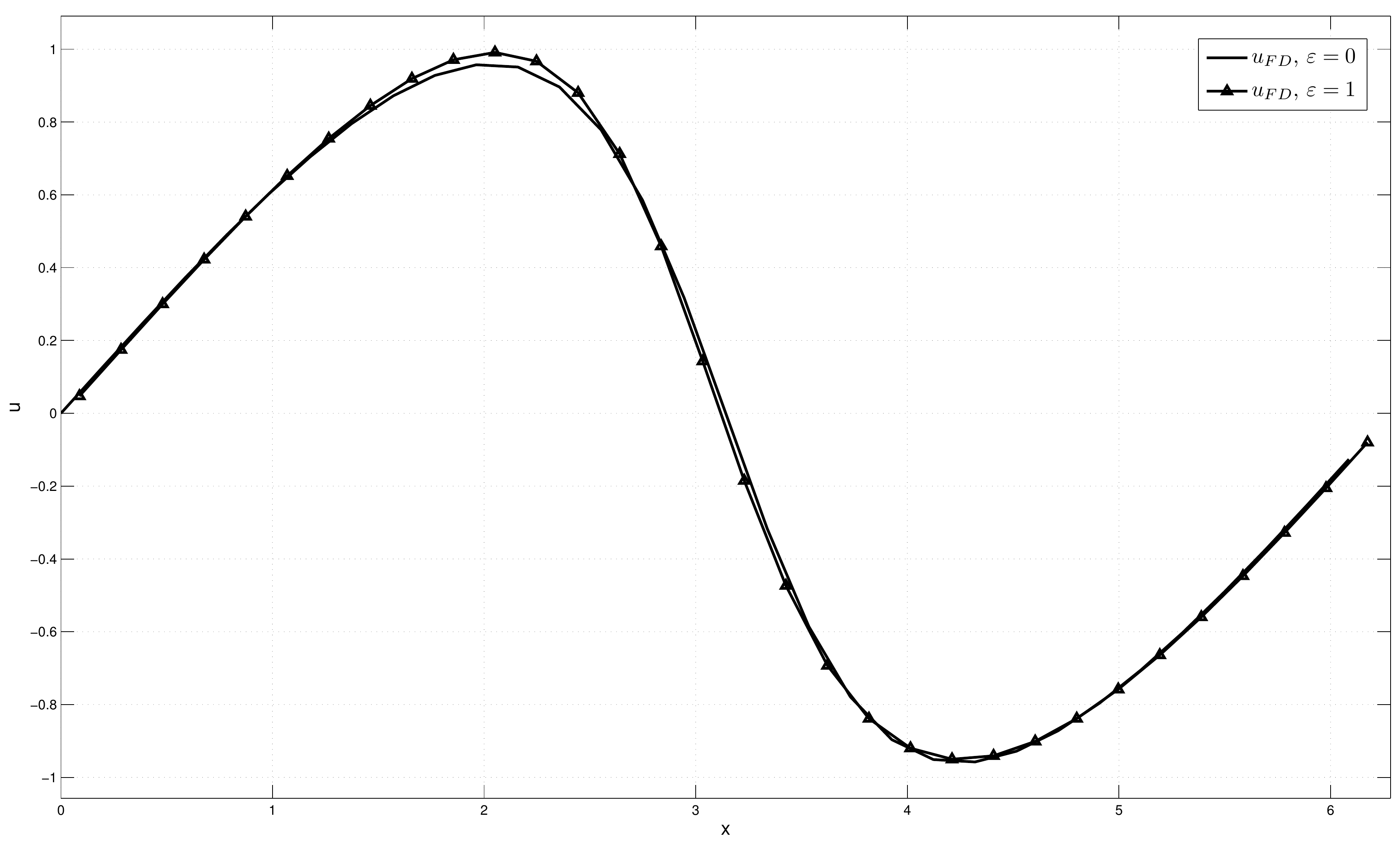}
 \caption{Integration using the classical FTCS discretization of Burgers equation~\eqref{eq:ViscousBurgers}. Solid line: Original integration in a resting reference frame. Solid lines with triangles: Integration in a reference frame moving with constant velocity $\ve=1$ as in~\cite{bern13Ay}. The results in the moving reference frames were shifted back to the origin for proper comparison.}
 \label{fig:NumericalVerificationGalileanInvarianceBurgers1}
\end{figure}

Unfortunately, to maintain Galilean invariance it is also not sufficient to carry out the numerical simulations with a standard finite difference scheme in a constantly moving reference frame as proposed in~\cite{bern13Ay}. It can be verified numerically that the resulting numerical solutions in the resting and in the convecting reference frames do not coincide, which is explicitly shown in Figure~\ref{fig:NumericalVerificationGalileanInvarianceBurgers1} for a FTCS discretization of Burgers equation. In this figure, we display the numerical solution at $t=0.5$ in the resting reference frame (solid line) and in a reference frame which moves with constant velocity $\ve_3=1$ (solid line with triangles) as in~\cite{bern13Ay}.

Instead of using a non-invariant finite difference scheme in a convecting reference frame, it is therefore desirable to construct proper finite difference discretizations that preserve the invariance group of a physical differential equations. The above observation on the incompatibility of stationary meshes with Galilean invariance have severe consequences on the design of finite difference models for the equations of fluid dynamics. In fact, it renders necessary to come up with strategies to combine the requirement of using moving meshes (in order to preserve Galilean invariance) with approaches that lead to discretization schemes having good numerical properties, such as stability, optimal grid adaptation (e.g.\ equidistribution of the discretization error) and the possibility for a parallel implementation. From a more general point of view, it is necessary to bridge the fields of group analysis and numerical analysis of differential equations.

To outline this connection for the example of Burgers equation considered in~\cite{bern13Ay} is the main aim of the present paper. In Section~\ref{sec:GalileanInvariantBurgers} we discuss invariant finite difference schemes for Burgers equation. We construct three different types of invariant numerical schemes, namely Lagrangian discretizations, invariant adaptive Eulerian schemes and invariant schemes employing an evolution--projection strategy. We relate these schemes to the remedy for reducing the effect of violation of Galilean invariance proposed in~\cite{bern13Ay}. Numerical results for the different schemes discussed are presented in Section~\ref{sec:NumericalResultsBurgers}. The final Section~\ref{sec:ConclusionBurgers} contains the conclusions of the paper.

\section{Invariant finite difference schemes for Burgers equations}\label{sec:GalileanInvariantBurgers}

As in~\cite{bern13Ay}, we introduce Burgers equation as a canonical example for high Reynolds number flows,\vspace{-0.5cm}
\begin{equation}\label{eq:ViscousBurgers}
 u_t+uu_x-\nu u_{xx}=0,
\end{equation}
where $\nu>0$ is the viscosity, which could be scaled to $1$ by means of an equivalence transformation. It is one of the most investigated models in the group analysis of differential equations, see e.g.~\cite{blum89Ay,blum10Ay,olve86Ay}. Its maximal Lie invariance algebra $\mathfrak g$ is spanned by the basis elements
\begin{equation}
 \p_t,\quad \p_x,\quad t\p_x+\p_u,\quad 2t\p_t+x\p_x-u\p_u,\quad t^2\p_t+tx\p_x+(x-tu)\p_u.
\end{equation}
The associated one-parameter Lie symmetry groups are
\begin{align}\label{eq:OneParameterGroupsViscousBurgers}
\begin{split}
 &\Gamma_1\colon\quad (t,x,u)\mapsto (t+\ve_1,x,u),\\
 &\Gamma_2\colon\quad (t,x,u)\mapsto (t,x+\ve_2,u),\\
 &\Gamma_3\colon\quad (t,x,u)\mapsto (t,x+\ve_3t,u+\ve_3),\\
 &\Gamma_4\colon\quad (t,x,u)\mapsto (e^{2\ve_4}t,e^{\ve_4}x,e^{-\ve_4}u),\\
 &\Gamma_5\colon\quad (t,x,u)\mapsto \left(\frac{t}{1-\ve_5t},\frac{x}{1-\ve_5t},u(1-\ve_5t)+\ve_5x\right),
\end{split}
\end{align}
showing that Burgers equation~\eqref{eq:ViscousBurgers} admits time translations, space translations, Galilean boosts, scalings and time inversions as one-parameter symmetry transformations.

Invariant numerical schemes for Eq.~\eqref{eq:ViscousBurgers} have already been investigated in the literature~\cite{chha10Ay,chha11Ay,doro11Ay,vali05Ay}. The schemes constructed in these references preserve the entire five-parameter symmetry group $G$ of Burgers equation. However, as was discussed in~\cite{bihl12By}, it is more natural to preserve only those symmetries that are compatible with a particular set of initial and boundary value problems chosen. In the present, we focus on periodic boundary conditions. The time inversion $\Gamma_5$ is then not compatible with a periodic domain as it does not map any periodic function $u$ to another periodic function for $\ve_5\ne0$. As a result we only aim to numerically preserve the first four symmetry transformations $\Gamma_1$--$\Gamma_4$. These transformations form the subgroup $G^1$ of the maximal Lie invariance group $G$ of Burgers equation. It is also important to note that the subgroup $G^1$ is typical for various models of fluid mechanics. Consequently, the strategies discussed below are also relevant for physically more interesting higher-dimensional models of hydrodynamics, such as the Euler or Navier--Stokes equations.

Arguably, the most important observation established in the field of invariant finite difference schemes is that it is generally not possible to maintain all symmetries of a system of differential equations if the discretization scheme is constructed on a fixed, orthogonal discretization mesh~\cite{doro11Ay,levi06Ay,rebe11Ay}. In the case of Burgers equation, it is the presence of the Galilean transformations $\Gamma_3$ that prohibits the use of a fixed discretization mesh. This was explicitly shown in the introduction. Hence, finite difference models operating on a fixed mesh cannot be Galilean invariant.

A possible remedy is to use the following expression as a discretization of Eq.~\eqref{eq:ViscousBurgers}
\begin{equation}\label{eq:InvariantSchemeBurgers}
 \frac{u^{n+1}_i-u^n_i}{\Delta t} + \left(u^n_i-\dot x^d_i\right)\frac{u^n_{i+1}-u^n_{i-1}}{x^n_{i+1}-x^n_{i-1}}-\frac{2\nu}{x^n_{i+1}-x^n_{i-1}} \left(\frac{u^n_{i+1}-u^n_i}{x^n_{i+1}-x^n_i}-\frac{u^n_{i}-u^n_{i-1}}{x^n_{i}-x^n_{i-1}}\right)=0
\end{equation}
where $\dot x^d_i=(x^{n+1}_i-x^n_i)/\Delta t $. Applying the transformations $\Gamma_1$--$\Gamma_4$ one readily verifies the invariance of this discretization. The reason for this scheme being invariant is that the introduced \textit{grid velocity} $\dot x^d_i$ transforms as $\dot x^d_i\to \dot x^d_i+\ve_4$ under the action of the Galilean transformation and the additional term involving $\ve_4$ is exactly compensated by the Galilean transformation of $u^n_i$. In other words, introducing a moving mesh into the discretization of Burgers equation restores the Galilean invariance of the finite difference scheme. See~\cite{doro11Ay,kim08Ay,olve01Ay,rebe11Ay} for further details on the systematic construction of invariant finite difference discretization schemes.

In~\cite{bihl12By} it was shown that the above discretization can be interpreted as a discretization of~\eqref{eq:ViscousBurgers} in terms of \textit{computational coordinates}, i.e.\ in a coordinate system that remains fixed in the presence of grid adaptation. To accomplish this transformation, one sets $x=x(\theta,\xi)$, where $\theta=t$ and $\xi=\xi(t,x)$ is the spatial computational coordinate. Transforming~\eqref{eq:ViscousBurgers} to the $(\theta,\xi)$ coordinates and discretizing the result using a FCTS scheme leads to~\eqref{eq:InvariantSchemeBurgers}.

The question remaining is how to determine the grid velocity $\dot x^d_i$, which involves the yet unknown location of the grid points on the subsequent time layer $n+1$. As pointed out in~\cite{rebe11Ay}, there are two main strategies to find $x^{n+1}_i$: The first is to construct a grid equation that is invariant under the same symmetry group as the discretization of the physical differential equation. The second method is to regard the grid adaptation as unconstrained from the symmetry requirements imposed by the physical differential equation, i.e.\ to use a non-invariant grid equation. We will mostly focus on the first method here as it is geometrically more grounded.

In the first method we require the grid equation to be invariant under the same symmetry group as is the discretization of the physical differential equation. In the present case, this amounts to constructing a grid equation that is invariant under the transformations $\Gamma_1$--$\Gamma_4$. One simple possibility is to take
\begin{equation}\label{eq:LagrangianGridEquationBurgers}
 \dot x^d_i-u_i^n=0,
\end{equation}
as this equation is obviously Galilean invariant and also does not violate the remaining transformations from $G^1$. This choice boils down to discretizing Burgers equation in \textit{Lagrangian coordinates}, i.e.\ the grid velocity equals the physical velocity. Indeed, Lagrangian discretization schemes are among the earliest examples of invariant discretizations for hydrodynamical equations, see e.g.~\cite{doro11Ay}. The problem with Lagrangian discretizations is that one in general does not have proper control over the evolution of the grid points. This can be a severe problem, especially in the multi-dimensional case, where grid points can concentrate in certain regions, deteriorating the local resolution of the scheme in areas away from these concentration regions.

Perhaps numerically more satisfying are grid equations that couple the evolution of the grid to the development of pronounced features in the numerical solution, i.e.\ to use proper grid adaptation strategies. Linking grid adaptation to the construction of invariant discretization schemes proved relevant in the numerical investigation of blow-up problems, see e.g.~\cite{budd96Ay,budd09Ay,huan10By}. A possible way to realize an \textit{invariant adaptive grid} is based on the equidistribution principle for a monitor function $\rho$,
\begin{equation}\label{eq:EquidistributionPrincipleBurgers}
 (\rho x_\xi)_\xi=0,
\end{equation}
which plays a central role in the construction of $r$-adaptive numerical schemes in one space dimension. See again~\cite{budd09Ay,huan10By} and references therein for an extensive discussion of the concept of equidistributing meshes. An invariant equidistributing mesh is obtained by discretizing~\eqref{eq:EquidistributionPrincipleBurgers} in a $G^1$-invariant way. The general feasibility of this approach depends on the structure of the symmetry group one aims to preserve~\cite{bihl12By} but in general it can be realized for the symmetry groups one usually encounters in hydrodynamics. The basis for this approach is to choose a proper monitor function $\rho$, that leads to a form of~\eqref{eq:EquidistributionPrincipleBurgers} that is invariant under the symmetry subgroup $G^1$ of Burgers equation and then to discretize this expression in an invariant way. A $G^1$-invariant monitor function is
\[
 \rho=\sqrt{1+\alpha u_x^2},
\]
which coincides with the arc-length function for $\alpha=1$. The reason for including a generic $\alpha$ in this expression is that the term $u_x$ is not scale invariant, i.e.\ it transforms as $\tilde u_{\tilde x}=e^{-2\ve_3}u_x$ and thus, as it stands, $\rho$ is not scale invariant. Extending the scalings of $(t,x,u)$ to an equivalence transformation involving $\alpha$ by adopting the transformation rule $\tilde \alpha=e^{2\ve_3}\alpha$ then indeed leads to a $G^1$-invariant function $\rho$. It should be stressed though that no such extension of the Galilean transformation to an equivalence transformation is needed to guarantee the Galilean invariance of the resulting form of the equidistribution principle~\eqref{eq:EquidistributionPrincipleBurgers}.

Using the modified arc-length weight function, a possible $G^1$-invariant discretization of~\eqref{eq:EquidistributionPrincipleBurgers} is
\begin{align}\label{eq:InvariantEquidistributionBurgers}
\begin{split}
 &(\rho^n_{i+1}+\rho^n_i)(x^{n+1}_{i+1}-x^{n+1}_i)-(\rho^n_i+\rho^n_{i-1})(x^{n+1}_i-x^{n+1}_{i-1})=0,\\ &\rho^n_i=\sqrt{1+\alpha\frac{u^n_{i+1}-u^n_{i-1}}{x^n_{i+1}-x^n_{i-1}}},
\end{split}
\end{align}
which can be solved using a relaxation scheme, such as e.g.\ Gau\ss-Seidel iteration to obtain $x^{n+1}_i$ and hence to complete the invariant numerical scheme~\eqref{eq:InvariantSchemeBurgers}. More sophisticated ways to solve~\eqref{eq:EquidistributionPrincipleBurgers} are conceivable as well and could be used to improve the quality of the resulting adaptive discretization scheme.

A further possibility for the construction of invariant numerical schemes is to invoke an \textit{evolution--projection strategy}. This idea was put forward in~\cite{nave10Ay,seib12Ay} for the non-invariant discretization of advection equations and extended in~\cite{bihl12Cy} to find an invariant discretization of the linear heat equation. The main approach in the evolution--projection strategy is to use the invariant numerical scheme and the invariant mesh equation only for a single integration step and to use a projection operator (i.e.\ an interpolation) to map the solution from the off-grid points back to the initial, uniformly spaced mesh. If the interpolation is done in an invariant way, i.e.\ the interpolation used preserves the invariance (sub)group of the system of differential equations being discretized, then the entire discretization procedure becomes invariant. The advantage of this approach is that moving meshes can be completely avoided.

In the present case of Burgers equation~\eqref{eq:ViscousBurgers} we observe that classical interpolation schemes such as linear, quadratic or cubic spline interpolation already preserve the invariance subgroup $G^1$. This means that we can use the aforementioned interpolations to re-map the solution $u^{n+1}_i$ defined at the points $x^{n+1}_i$ back to $\hat x_i^{n+1}\in \{x_i^n\}$, without breaking the invariance of the scheme.

We show this explicitly for linear interpolation here, which is defined as
\[
 u(\hat x^{n+1}_i)=u_i^{n+1}+\frac{u_{i+1}^{n+1}-u_i^{n+1}}{x_{i+1}^{n+1}-x_i^{n+1}}(\hat x^{n+1}_i-x_i^{n+1})=\mathcal L(x_i^{n+1},x_{i+1}^{n+1},u_i^{n+1},u_{i+1}^{n+1};\hat x_i^{n+1})
\]
for the interpolation of values $\hat x_i^{n+1}$ lying within the interval $[x_i^{n+1},x_{i+1}^{n+1}]$. Then, for transformations of the form $\Gamma_1$--$\Gamma_4$, we obtain that
\[
 \tilde u(\widetilde{\hat x}^{n+1}_i)-\mathcal L(\tilde x_i^{n+1},\tilde x_{i+1}^{n+1},\tilde u_i^{n+1},\tilde u_{i+1}^{n+1};\widetilde{\hat x}_i^{n+1})= u(\hat x^{n+1}_i)-\mathcal L(x_i^{n+1},x_{i+1}^{n+1},u_i^{n+1},u_{i+1}^{n+1};\hat x_i^{n+1}),
\]
which proves the invariance of linear interpolation under transformations from $G^1$. Similarly, invariance of quadratic and cubic spline interpolation can be shown.

Regarding the use of non-invariant grid equations, which is the second possibility to complete the description of the scheme~\eqref{eq:InvariantSchemeBurgers}, in principle all choices excluding $x^{n+1}_i=x^n_i$ are admissible. We will now return to the remedy proposed in~\cite{bern13Ay}. The proposed approach falls into the category of non-invariant grid equations. In that paper, it was suggested to use a reference frame moving with the (constant) bulk velocity of the flow. In the case of Burgers equation, the authors set
\begin{equation}\label{eq:BernardiGridEquation}
 x^{n+1}_i=x^n_i+c\Delta t,
\end{equation}
where $c=\const$. It is readily verified that this grid equation is \textit{not} Galilean invariant, as
\[
 \tilde x^{n+1}_i-\tilde x^n_i-c\Delta \tilde t= x^{n+1}_i-x^n_i-(c-\ve_4)\Delta t.
\]
If $c\ne0$ is the bulk velocity of a flow, one may indeed expect that the numerical results obtained from scheme~\eqref{eq:InvariantSchemeBurgers} with grid equation~\eqref{eq:BernardiGridEquation} are better than for scheme~\eqref{eq:InvariantSchemeBurgers} on a stationary grid. More precisely, \textit{Galilean invariance could be restored by extending the transformation $\Gamma_4$ to $c$ by setting $\tilde c=c+\ve_4$.} This extension is justified in case $c$ is related to $u$, which is the main reason why the remedy proposed in~\cite{bern13Ay} may work. In other words, reference frames moving with a constant velocity could be made Galilean invariant in the sense that Galilean transformations have an extension to equivalence transformations for such reference frames. A similar extension of the transformation $\Gamma_3$ to $c$ is necessary to incorporate the scale invariance.

For the sake of convenience, we summarized the characteristics of the different schemes discussed in the present section in Table~\ref{tab:OverviewSchemesBurgers}.

\begin{table}[!ht]
\centering
\caption{Different numerical schemes for the viscous Burgers equation~\eqref{eq:ViscousBurgers}.}
\begin{tabular}{|c|c|c|c|}
\hline
 Numerical scheme & Grid equation & Galilean invariance & Grid spacing \\
\hline
Finite differences & None & Not invariant & constant\\
Lagrangian & Lagrangian grid movement & fully invariant & variable\\
Eulerian adaptive & Equidistribution principle & fully invariant & variable\\
\cite{bern13Ay} Bernardi et al.& Constant grid movement  & semi-invariant & constant\\
Evolution--projection & Lagrangian grid movement  & fully invariant & constant\\
\hline
\end{tabular}
\label{tab:OverviewSchemesBurgers}
\end{table}

\section{Numerical results}\label{sec:NumericalResultsBurgers}

In this section we present some numerical results obtained from the invariant discretization schemes introduced in the previous section. For all the experiments, we use $u(0,x)=\sin(x)$ as the initial condition on a $2\pi$-periodic domain and the integration is carried out up to time $t=0.5$ for $\nu=0.1$. Unless otherwise stated, we use $N=64$ grid points and fix the time step with $\Delta t \propto h^2$, where $h$ is the mean grid spacing over the domain.

\begin{figure}[!ht]
 \centering
 \includegraphics[scale=0.40]{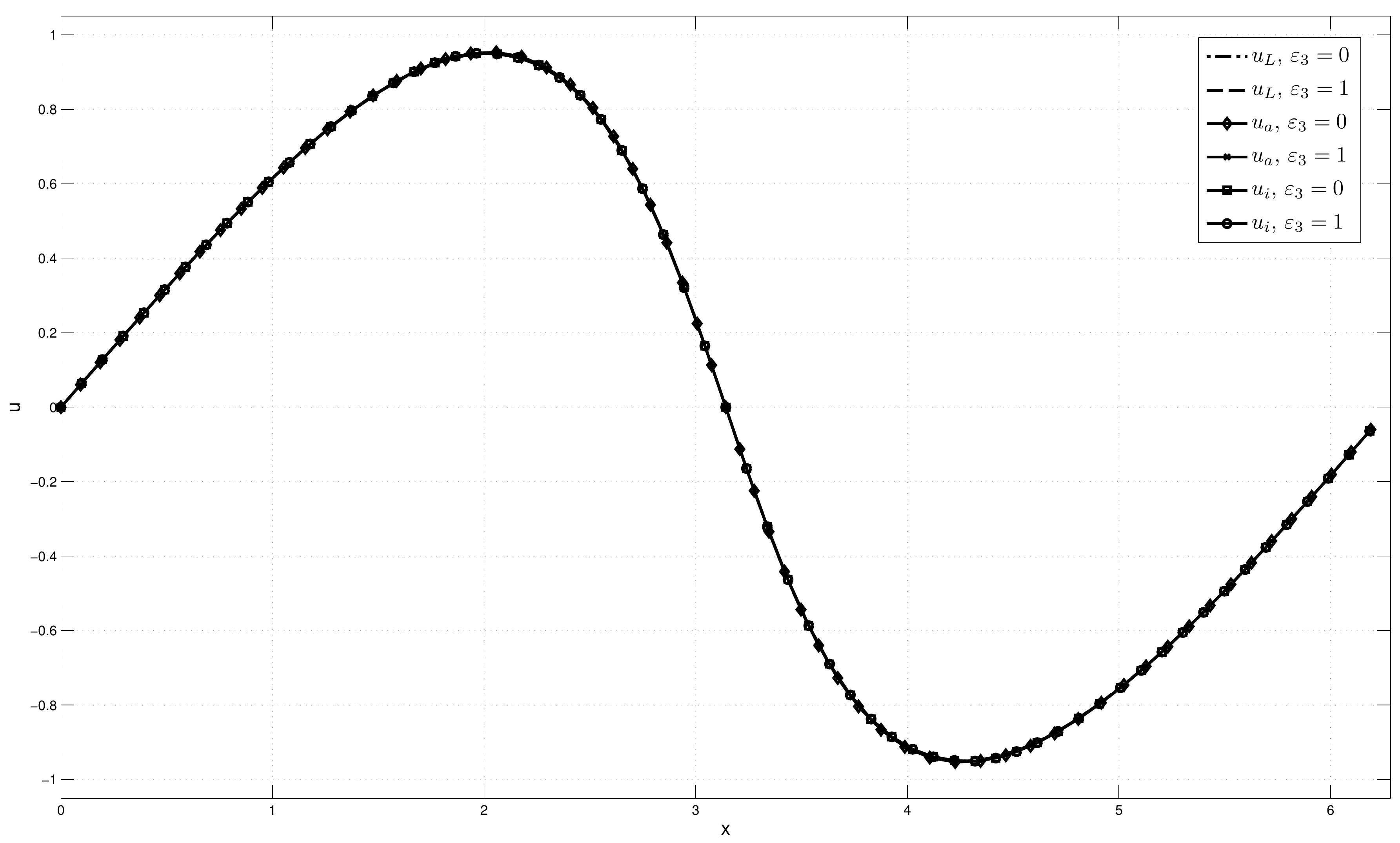}
 \caption{Integration of Burgers equation~\eqref{eq:ViscousBurgers} using the three invariant discretization schemes introduced in Section~\ref{sec:GalileanInvariantBurgers}. Dashed--dotted line: Lagrangian scheme in the resting reference frame. Dashed line: Lagrangian scheme in a reference frame moving with constant velocity $\ve_3=1$. Solid line with diamonds: Eulerian adaptive scheme in the resting reference frame. Solid line with crosses: Eulerian adaptive scheme in a reference frame moving with constant velocity $\ve_3=1$. Solid line with squares: Evolution--projection scheme in the resting reference frame. Solid line with circles: Evolution--projection scheme in a reference frame moving with constant velocity $\ve_3=1$.}
 \label{fig:NumericalVerificationGalileanInvarianceBurgers2}
\end{figure}

In Figure~\ref{fig:NumericalVerificationGalileanInvarianceBurgers2} we carry out numerical integrations for three invariant numerical schemes for Burgers equation based on~\eqref{eq:InvariantSchemeBurgers} and employing different grid equations: (i) The fully Lagrangian scheme uses the Lagrangian grid equation~\eqref{eq:LagrangianGridEquationBurgers}, (ii) in the invariant Eulerian adaptive scheme the grid points $\{x_i^{n+1}\}$ is determined from the invariant discretization of the equidistribution principle~\eqref{eq:EquidistributionPrincipleBurgers} and (iii) in the evolution--projection scheme we use the Lagrangian grid equation~\eqref{eq:LagrangianGridEquationBurgers} and use quadratic interpolation to re-map the off-grid points to their original location at the previous time step.

For all schemes we numerically verify Galilean invariance, i.e.\ each of the pairs of integrations in a resting and a constantly moving coordinate system yields visually the same numerical solution for (i) the Lagrangian scheme (dashed--dotted line, $\ve_3=0$ and dashed line, $\ve_3=1$), (ii) the Eulerian adaptive scheme (solid line with diamonds, $\ve_3=0$ and solid line with crosses $\ve_3=1$) and (iii) the evolution--projection scheme (solid line with squares, $\ve_3=0$ and solid line with circles $\ve_3=1$). Moreover, it is seen from Figure~\ref{fig:NumericalVerificationGalileanInvarianceBurgers2} that all three schemes yield approximately the same numerical solution. We aim to detail and analyze this result further now.

\begin{figure}[!ht]
 \centering
 \includegraphics[scale=0.4]{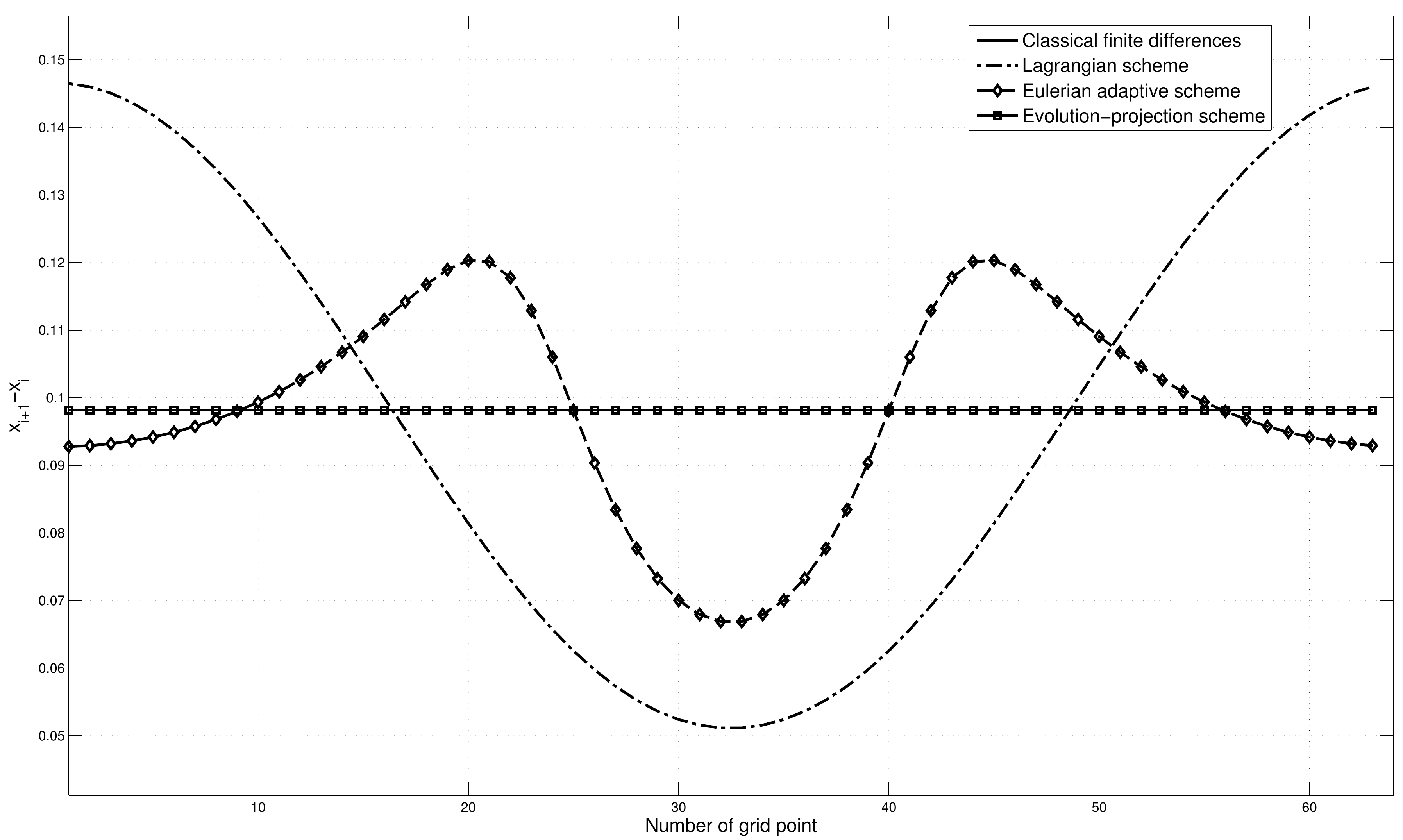}
 \caption{Grid spacing $\Delta x^n_i=x^n_{i+1}-x^n_i$ at final time $t=0.5$ for the four discretizations of Burgers equation shown in Figures~\ref{fig:NumericalVerificationGalileanInvarianceBurgers1} and \ref{fig:NumericalVerificationGalileanInvarianceBurgers2}. }
 \label{fig:NumericalVerificationGalileanInvarianceBurgers3}
\end{figure}

The difference between the three invariant numerical schemes for Burgers equation is the invoked grid equation. For the sake of a clearer presentation, we depict the grid spacing $\Delta x^n_i=x^n_{i+1}-x^n_i$ as a function of the location of the grid points $x^n_i$ at the final integration time $t=0.5$ in Figure~\ref{fig:NumericalVerificationGalileanInvarianceBurgers3}. For the classical, non-invariant finite difference scheme the spacing is by definition constant (solid line). For the Lagrangian discretization (dashed line) the location of the grid points depends on the solution itself and therefore one does not have control over the local resolution. It is a mere consequence of equating the physical velocity and the grid velocity. In the adaptive Eulerian scheme (dashed line with diamonds) we observe a proper concentration of the grid points along the building shock, which follows from using the equidistribution principle and the arc-length monitor function. Away from the steepening front, the points remain quasi-equally distributed. By construction, the grid points in the evolution--projection method (solid line with squares) are again equally spaced.

\begin{figure}[!ht]
 \centering
 \includegraphics[scale=0.4]{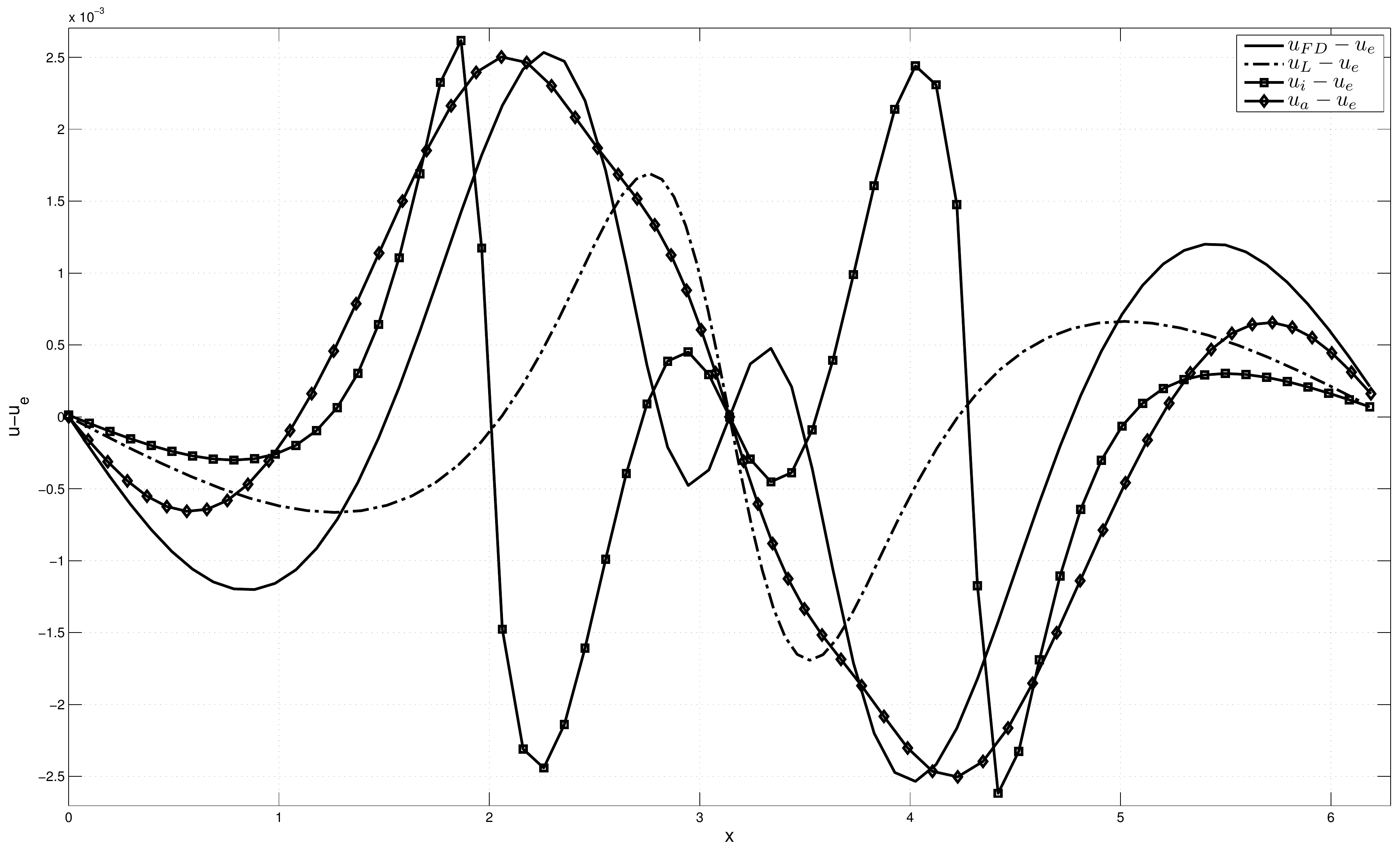}
 \caption{Pointwise difference of the numerical and the exact solution at final time $t=0.5$ for the four discretizations of Burgers equation shown in Figures~\ref{fig:NumericalVerificationGalileanInvarianceBurgers1} and \ref{fig:NumericalVerificationGalileanInvarianceBurgers2}.}
 \label{fig:NumericalVerificationGalileanInvarianceBurgers4}
\end{figure}

To estimate the overall accuracy of the different schemes, in Figure~\ref{fig:NumericalVerificationGalileanInvarianceBurgers4} we display the pointwise differences of the numerical solutions and the exact solution of the chosen initial value problem $u(0,x)=\sin(x)$ for Burgers equation, which is
\begin{align*}
 & u_e(t,x)=2\nu\dfrac{\sum_{j=1}^\infty a_jje^{-\nu tj^2}\sin jx}{\sum_{j=0}^\infty a_je^{-\nu tj^2}\cos jx},\\
 &a_0=\frac{1}{2\pi}\int_{0}^{2\pi}e^{-(1-\cos x)/(2\nu)}\ddd x,\quad a_{j>0}=\frac{1}{\pi}\int_{0}^{2\pi}e^{-(1-\cos x)/(2\nu)}\cos j x\, \ddd x.
\end{align*}
It can be seen from Figure~\ref{fig:NumericalVerificationGalileanInvarianceBurgers4} that all the numerical schemes achieve a comparable accuracy. The overall $l_\infty$-errors of the schemes for the runs depicted in Figure~\ref{fig:NumericalVerificationGalileanInvarianceBurgers4} are summarized in Table~\ref{tab:RMSerrorsBurgers}.

\begin{table}[!ht]
\centering
\vspace{-0.4cm}
\caption{$l_\infty$-errors for the various discretization schemes for the viscous Burgers equation~\eqref{eq:ViscousBurgers} with $N=64$ grid points.}
\begin{tabular}{|c|c|c|c|c|}
\hline
 & classical FD & Lagrangian & Eulerian adaptive & Evolution--projection \\
\hline
$||E||_{l_{\infty}}$ & $2.53\cdot 10^{-3}$ & $1.69\cdot 10^{-3}$ & $2.50\cdot 10^{-3}$ & $2.63\cdot 10^{-3}$\\
\hline
\end{tabular}
\label{tab:RMSerrorsBurgers}
\end{table}

\begin{figure}[!ht]
 \centering
 \includegraphics[scale=0.4]{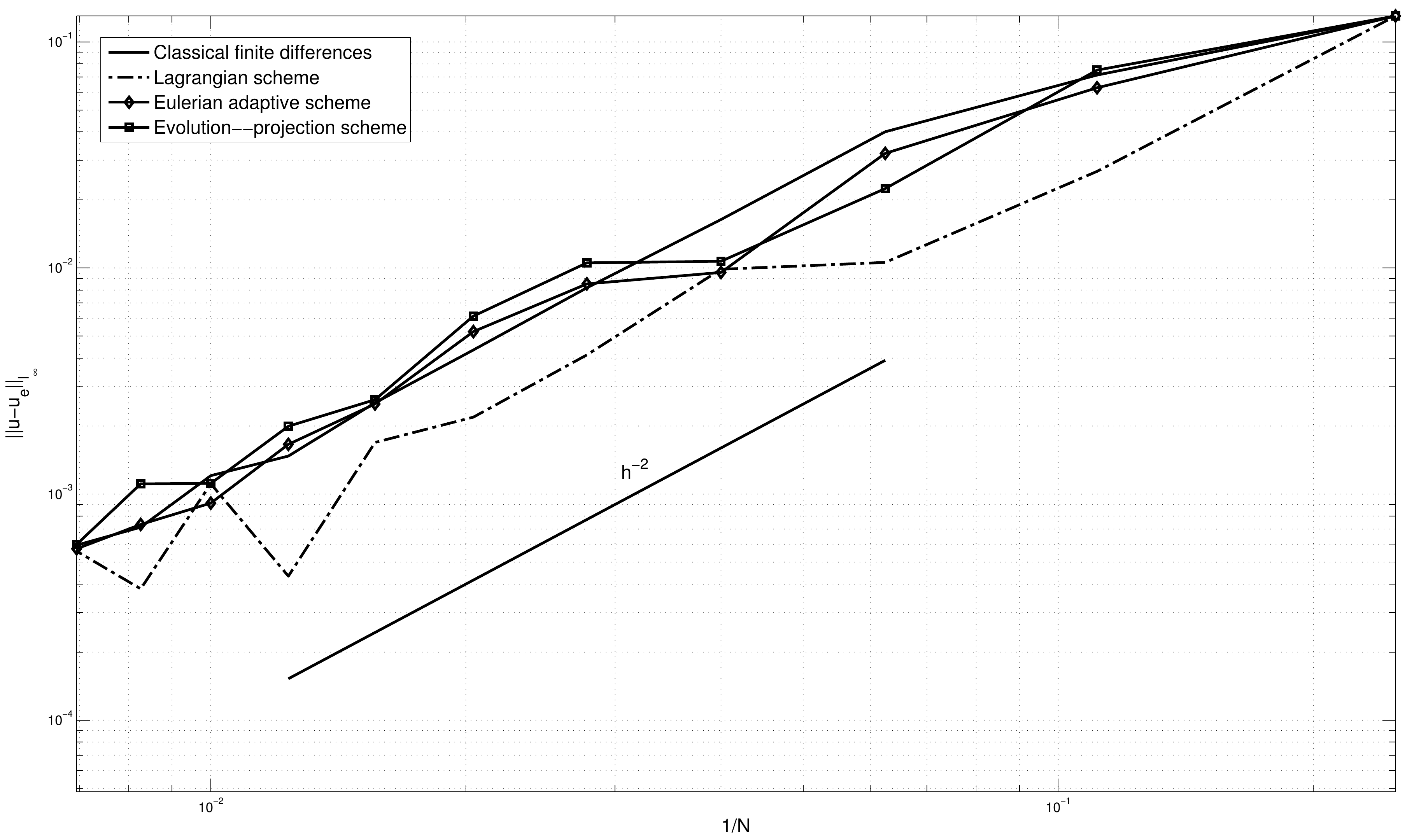}
 \caption{Convergence plots for the four numerical schemes presented above for $N\in\{4,8,16,32,64,128,256,512\}$ grid points. Solid line: Classical finite difference scheme. Dashed--dotted line: Lagrangian scheme. Solid line with diamonds: Eulerian adaptive scheme. Solid line with squares: Evolution--projection scheme.}
 \label{fig:NumericalVerificationGalileanInvarianceBurgers5}
\end{figure}

The convergence rates of the different schemes in the $l_\infty$-norm are depicted in Figure~\ref{fig:NumericalVerificationGalileanInvarianceBurgers5} using $N\in\{4,8,16,32,64,128,256,512\}$ grid points. It is seen from this plot that all schemes yield comparable errors. The overall convergence rates demonstrate that all three types of invariant schemes for Burgers equation introduced in this paper are asymptotically of second order and that the different invariant discretization strategies do not alter the accuracy of the underlying approximations. More details on this can be found in~\cite{bihl12Cy}.

\section{Conclusion}\label{sec:ConclusionBurgers}

In the present paper we have revisited the problem recently pointed out in~\cite{bern13Ay} that classical finite difference discretizations of the governing equations of hydrodynamics violate Galilean invariance. We have discussed three possible ways of constructing Galilean invariant finite difference schemes for Burgers equation: (i) Lagrangian discretization schemes, (ii) Eulerian adaptive discretizations and (iii) discretizations using an evolution--projection strategy. These three approaches can be readily adapted to the two- and three-dimensional Euler or Navier--Stokes equations. The approach proposed in~\cite{bern13Ay} leads to a semi-invariant scheme in that the authors use the discretization on a moving coordinate system but employ a non-invariant grid equation by moving the grid points with a constant velocity (e.g.\ the bulk velocity) rather than with the actual physical velocity. If properly done, this approach can indeed reduce the effects of the violation of Galilean invariance in classical finite difference discretizations, while still not a fully Galilean invariant schemes.

All of the invariant discretization approaches presented above have their advantages and disadvantages. Purely Lagrangian discretization schemes are not in widespread use as such schemes usually lead to a strong concentration of grid points, leaving other regions of the domain poorly resolved. Moreover, in multi-dimensional cases of interest in hydrodynamics, Lagrangian schemes can lead to tangled meshes. At the same time, Lagrangian schemes are good in that they are able to preserve sharp interfaces within a fluid. The Eulerian adaptive approach is attractive because it provides a natural way to link the problem of finding invariant discretization schemes to the properties of the numerical solution at each time step. At the same time, the computational overhead required to efficiently generate the meshes can be a crucial factor determining the feasibility of the adaptation methodology, especially for multi-dimensional systems. Finally, in situations where adaptive numerical schemes are not desirable, the evolution--projection strategy based on the Lagrangian grid equation (or any other invariant grid equation one is able to find) and a simple (but invariant) interpolation is a possible way to maintain Galilean invariance in a finite difference scheme while still being able to operate on a fixed, uniformly spaced mesh. The drawback of the evolution--projection approach is that it also requires an additional operation, the interpolation, which may cause additional computational overhead compared to, for instance, the Lagrangian scheme.

Applied to Burgers equation, the three invariant discretization methodologies yielded numerical schemes that are asymptotically second order accurate. Compared to the standard FTCS scheme, which is also second order accurate, one thus gets as a bonus the preservation of a geometric feature of Burgers equation, namely a subgroup of its symmetry group. The preservation of this symmetry subgroup can be a crucial factor for problems for which it is inevitable to carry out the simulations in a moving reference frame. The construction of higher order invariant numerical schemes is the subject of current research and will be reported in the future.

\section*{Acknowledgements}

This research was supported by the Austrian Science Fund (FWF), project J3182--N13 (AB). JCN wishes to acknowledge partial support from the NSERC Discovery Program, and the National Science Foundation through grant DMS-0813648.

{\footnotesize\setlength{\itemsep}{0ex}

}


\begin{thebibliography}{10}
\providecommand{\url}[1]{\texttt{#1}}
\providecommand{\urlprefix}{URL }
\expandafter\ifx\csname urlstyle\endcsname\relax
  \providecommand{\doi}[1]{doi:\discretionary{}{}{}#1}\else
  \providecommand{\doi}{doi:\discretionary{}{}{}\begingroup
  \urlstyle{rm}\Url}\fi
\providecommand{\eprint}[2][]{\url{#2}}

\bibitem{bern13Ay}
Bernardini M., Pirozzoli S., Quadrio M. and Orlandi P., Turbulent channel flow
  simulations in convecting reference frames, \emph{J. Comput. Phys.}
  \textbf{232} (2013), 1--6.

\bibitem{bihl12Cy}
Bihlo A. and Nave J.-C., Invariant discretization schemes for the heat equation,
  2012, arXiv:1209.5028, 21 pp.

\bibitem{bihl12By}
Bihlo A. and Popovych R.O., Invariant discretization schemes for the
  shallow-water equations, \emph{SIAM J. Sci. Comput.} \textbf{34} (2012),
  B810--B839, arXiv:1201.0498.

\bibitem{blum89Ay}
Bluman G. and Kumei S., \emph{{Symmetries and differential equations}},
  Springer, New York, 1989.

\bibitem{blum10Ay}
Bluman G.W., Cheviakov A.F. and Anco S.C., \emph{Application of symmetry
  methods to partial differential equations}, Springer, New York, 2010.

\bibitem{budd01Ay}
Budd C. and Dorodnitsyn V.A., {Symmetry-adapted moving mesh schemes for the
  nonlinear Schr\"{o}dinger equation}, \emph{J. Phys. A} \textbf{34} (2001),
  10387--10400.

\bibitem{budd96Ay}
Budd C.J., Huang W. and Russell R.D., Moving mesh methods for problems with
  blow-up, \emph{SIAM J. Sci. Comput.} \textbf{17} (1996), 305--327.

\bibitem{budd09Ay}
Budd C.J., Huang W. and Russell R.D., Adaptivity with moving grids, \emph{Acta
  Numer.} \textbf{18} (2009), 111--241.

\bibitem{chha10Ay}
Chhay M. and Hamdouni A., A new construction for invariant numerical schemes
  using moving frames, \emph{C.~R. Mecanique} \textbf{338} (2010), 97--101.

\bibitem{chha11Ay}
Chhay M., Hoarau E., Hamdouni A. and Sagaut P., Comparison of some
  {L}ie-symmetry-based integrators, \emph{J. Comput. Phys.} \textbf{230}
  (2011), 2174--2188.

\bibitem{doro11Ay}
Dorodnitsyn V., \emph{Applications of Lie Groups to Difference Equations},
  vol.~8 of \emph{Differential and integral equations and their applications},
  Chapman \& Hall/CRC, Boca Raton, FL, 2011.

\bibitem{doro03Ay}
Dorodnitsyn V.A. and Kozlov R., A heat transfer with a source: the complete set
  of invariant difference schemes, \emph{J. Nonlin. Math. Phys.} \textbf{10}
  (2003), 16--50.

\bibitem{huan10By}
Huang W. and Russell R.D., \emph{Adaptive Moving Mesh Methods}, Springer, New
  York, 2010.

\bibitem{kim08Ay}
Kim P., {Invariantization of the Crank--Nicolson method for Burgers' equation},
  \emph{Physica D} \textbf{237} (2008), 243--254.

\bibitem{levi06Ay}
Levi D. and Winternitz P., {Continuous symmetries of difference equations},
  \emph{J. Phys. A} \textbf{39} (2006), R1--R63.

\bibitem{nave10Ay}
Nave J.-C., Rosales R.R. and Seibold B., A gradient-augmented level set method
  with an optimally local, coherent advection scheme, \emph{J. Comput. Phys.}
  \textbf{229} (2010), 3802--3827.

\bibitem{olve86Ay}
Olver P.J., \emph{{Application of Lie groups to differential equations}},
  Springer, New York, 2000.

\bibitem{olve01Ay}
Olver P.J., Geometric foundations of numerical algorithms and symmetry,
  \emph{Appl. Algebra Engrg. Comm. Comput.} \textbf{11} (2001), 417--436.

\bibitem{rebe11Ay}
Rebelo R. and Valiquette F., Symmetry preserving numerical schemes for partial
  differential equations and their numerical tests, {\it J. Difference Equ.
  Appl.}, to appear, arXiv:1110.5921, 2012.

\bibitem{seib12Ay}
Seibold B., Rosales R.R. and Nave J.-C., Jet schemes for advection problems,
  \emph{Discrete Contin. Dyn. Syst. Ser. B} \textbf{17} (2012), 1229--1259.

\bibitem{vali05Ay}
Valiquette F. and Winternitz P., {Discretization of partial differential
  equations preserving their physical symmetries}, \emph{J. Phys. A}
  \textbf{38} (2005), 9765--9783.

\end{thebibliography}
\end{document}